\documentclass{amsart}

\usepackage{amsmath,amssymb,amsthm,mathtools}
\usepackage{bm}
\usepackage{enumitem}
\usepackage[expansion=false]{microtype}
\usepackage[colorlinks=true,linkcolor=red,citecolor=blue,urlcolor=blue]{hyperref}
\usepackage[margin=0.8in]{geometry}

\newtheorem{theorem}{Theorem}[section]
\newtheorem{proposition}[theorem]{Proposition}
\newtheorem{corollary}[theorem]{Corollary}
\newtheorem{lemma}[theorem]{Lemma}

\newtheorem{example}[theorem]{Example}

\newcommand{\D}{\mathcal P}
\newcommand{\M}{\mathbb{M}_{n}}

\title[Power Moduli and Sharp Hybrid Matrix Triangle Inequalities]{Power Moduli and Sharp Hybrid Matrix Triangle Inequalities}

\author{Mohammad Alakhrass}
\address{Department of Mathematics, University of Sharjah, Sharjah 27272, UAE}
\email{malakhrass@sharjah.ac.ae}

\subjclass[2020]{15A42, 15A60, 47A30, 47A60.}
\keywords{Power moduli, symmetric modules, Olson's spectral order, hybrid matrix triangle inequalities, matrix geometric mean, weak log-majorization, binormal matrices}
\begin{document}

\begin{abstract}
We prove sharp hybrid triangle inequalities in the Loewner order
for sums of complex matrices in terms of their power moduli. The
bounds involve the matrix geometric mean and a unitary polar
factor of the sum. As consequences, we obtain singular value
inequalities, weak log-majorization, and bounds for all unitarily
invariant norms. The arithmetic, quadratic, and maximal moduli
arise as special cases. When the summands are binormal, the bounds
are nonincreasing in the exponent and can improve strictly on the
quadratic bound. An example shows that this monotonicity can fail
for general summands.
\end{abstract}

\maketitle

\section{Introduction}

Let \(M_n(\mathbb C)\) denote the algebra of complex \(n\times n\)
matrices. For Hermitian \(X,Y\in M_n(\mathbb C)\), we write
\(X\le Y\) when \(Y-X\) is positive semidefinite. For
\(A\in M_n(\mathbb C)\), write
\[
|A|=(A^*A)^{1/2},
\qquad
|A^*|=(AA^*)^{1/2}.
\]
For \(p>0\), define the power modulus of \(A\) by
\begin{equation}\label{definitions}
|A|_p
=
\left(\frac{|A|^p+|A^*|^p}{2}\right)^{1/p}.
\end{equation}
This is the power mean of the two positive semidefinite matrices
\(|A|\) and \(|A^*|\). See \cite{Bhatia-2012} for more details about the power mean for any two positive semidefinite matrices. At \(p=1\) and
\(p=2\), it gives the symmetric and quadratic symmetric moduli,
respectively:
\[
|A|_1=|A|_{\mathrm{sym}}
=\frac{|A|+|A^*|}{2},
\qquad
|A|_2=|A|_{\mathrm{qsym}}
=\left(\frac{|A|^2+|A^*|^2}{2}\right)^{1/2}.
\]
The limit as \(p\to\infty\) exists in norm and gives the maximal
modulus associated with Olson's spectral order:
\[
|A|_\infty
:=\lim_{p\to\infty}|A|_p
=\lim_{p\to\infty}\bigl(|A|^p+|A^*|^p\bigr)^{1/p}
=|A|_\vee.
\]
For invertible \(A\), the formula in \eqref{definitions} also
defines \(|A|_p\) for \(p<0\), and its continuous extension at
\(p=0\) is
\begin{equation}\label{zero}
|A|_0
:=\lim_{p\to0}|A|_p
=\exp\!\left(\frac{\log|A|+\log|A^*|}{2}\right).
\end{equation}

By Theorem~\ref{N-Mo}, the power modulus is nondecreasing
in the Loewner order on \([1,\infty]\). More precisely,

\begin{equation}\label{intro:power-monotonicity}
|A|_p\le |A|_q,
\qquad 1\le p\le q\le\infty.
\end{equation}
In particular, every power modulus with \(1\le p\le2\) lies
between \(|A|_{\mathrm{sym}}\) and \(|A|_{\mathrm{qsym}}\).

Let \(\|\cdot\|\) denote an arbitrary unitarily invariant norm.
Zhang \cite{Zhang-2026} established the comparisons
\begin{align}
\tfrac12\|A\|
&\le \bigl\||A|_{\mathrm{sym}}\bigr\|
\le \|A\|,
\label{Zahng 1}\\
2^{-1/2}\|A\|
&\le \bigl\||A|_{\mathrm{qsym}}\bigr\|
\le 2^{1/2}\|A\|,
\label{Zahng 2}\\
\bigl\||A|_{\mathrm{sym}}\bigr\|
&\le \bigl\||A|_{\mathrm{qsym}}\bigr\|
\le 2^{1/2}\bigl\||A|_{\mathrm{sym}}\bigr\|.
\label{Zahng 3}
\end{align}
Zhang also proved the triangle inequalities
\begin{align}
\left\|\sum_{k=1}^m A_k\right\|
&\le 2\left\|\sum_{k=1}^m|A_k|_{\mathrm{sym}}\right\|,
\label{Zahng 4}\\
\left\|\sum_{k=1}^m A_k\right\|
&\le \sqrt{2}
\left\|\sum_{k=1}^m|A_k|_{\mathrm{qsym}}\right\|.
\label{Zahng 5}
\end{align}
The constants are sharp uniformly over all unitarily invariant
norms. Corollary~\ref{cor:power-norm-chains} gives the following
comparisons with the power moduli:
\begin{align*}
\bigl\||A|_{\mathrm{sym}}\bigr\|
&\le \bigl\||A|_p\bigr\|
\le \bigl\||A|_{\mathrm{qsym}}\bigr\|
\le \sqrt{2}\,\bigl\||A|_{\mathrm{sym}}\bigr\|,
&&1\le p\le2,\\
\bigl\||A|_{\mathrm{sym}}\bigr\|
&\le \bigl\||A|_{\mathrm{qsym}}\bigr\|
\le \bigl\||A|_p\bigr\|
\le 2^{1-1/p}\bigl\||A|_{\mathrm{sym}}\bigr\|,
&&2\le p\le\infty.
\end{align*}

Bourin and Lee \cite{BL-hybrid} refined \eqref{Zahng 5} in the
Loewner order. For some unitary \(V_2\),
\begin{equation}\label{Bourn 1}
\left|\sum_{k=1}^m A_k\right|\le\sqrt{2}\biggl(\Bigl(\sum_{k=1}^m |A_k|_{2}\Bigr)\#\Bigl(\sum_{k=1}^m V_2^*|A_k|_{2}V_2\Bigr)\biggr).
\end{equation}
Here \(\#\) denotes the matrix geometric mean. They also proved
the maximal-modulus inequality
\begin{equation}\label{Bourn 2}
\left|\sum_{k=1}^m A_k\right|\le \biggl(\Bigl(\sum_{k=1}^m |A_k|_{\infty}\Bigr)\#\Bigl(\sum_{k=1}^m V_\infty^*|A_k|_{\infty}V_\infty\Bigr)\biggr)
\end{equation}
for a unitary \(V_\infty\).

We extend these bounds to \(1\le p\le\infty\). If \(V\) is a
unitary polar factor of \(\sum_{k=1}^m A_k\), then
\begin{equation}\label{intro-power-triangle}
\begin{aligned}
\left|\sum_{k=1}^m A_k\right|&\le 2^{1/p}\biggl(\Bigl(\sum_{k=1}^m |A_k|_{p}\Bigr)\#\Bigl(\sum_{k=1}^m V^*|A_k|_{p}V\Bigr)\biggr)\\
&\le \frac{2^{1/p}}{2}\sum_{k=1}^m\bigl(|A_k|_p+V^*|A_k|_pV\bigr).
\end{aligned}
\end{equation}
We use the convention \(2^{1/\infty}=1\). The coefficient \(2^{1/p}\) is sharp.
At \(p=1\), \eqref{intro-power-triangle} gives a Loewner-order
refinement of \eqref{Zahng 4}. At \(p=2\) and \(p=\infty\), it
recovers \eqref{Bourn 1} and \eqref{Bourn 2}, respectively.

The geometric-mean bound yields
\[
\left|\sum_{k=1}^m A_k\right|\prec_{w\log}2^{1/p}\sum_{k=1}^m |A_k|_{p},
\qquad
\left\|\sum_{k=1}^m A_k\right\|\le 2^{1/p}\Bigl\|\sum_{k=1}^m |A_k|_{p}\Bigr\|,
\]
where \(\prec_{w\log}\) denotes weak log-majorization. For \(m=1\),
the cases \(p=1\) and \(p=2\) refine the lower norm bounds in
\eqref{Zahng 1} and \eqref{Zahng 2}. At \(p=\infty\), we have
\[
|A|\le |A|_\vee,\qquad \|A\|\le\bigl\||A|_\vee\bigr\|.
\]

We show that \(2^{1/p}\sum_{k=1}^m |A_k|_p\) is nonincreasing
in the Loewner order when each summand \(A_k\) is binormal. In this case,
larger exponents can improve the quadratic bound. An example
with two matrices shows that the operator and
norm bounds can increase from \(p=2\) to \(p=4\) for general
summands.

Our proof uses the positive-block method employed by Bourin and
Lee \cite{BL-hybrid}. The same method was used in
\cite{Alakhrass02021} to derive inequalities for sectorial matrices.
Lemma~\ref{lem:positive-blocks} shows how positivity of a block
matrix gives the geometric-mean bound and its singular value
consequences.

\section{Positive block matrices and the geometric mean}

For a Hermitian matrix \(X\), we order its eigenvalues as
\(\lambda_1(X)\ge\cdots\ge\lambda_n(X)\). For an arbitrary matrix
\(X\), let \(s_j(X)=\lambda_j(|X|)\), and write
\(s(X)=(s_1(X),\ldots,s_n(X))\). We use \(\|\cdot\|\) for an
arbitrary unitarily invariant norm.

For \(X>0\) and \(Y >0\), the matrix geometric mean is given by
\[
X\#Y
=
X^{1/2}\bigl(X^{-1/2}YX^{-1/2}\bigr)^{1/2}X^{1/2}.
\]
For general \(X,Y\ge0\), it is defined by the limit
\[
X\#Y=\lim_{\varepsilon\downarrow0}
(X+\varepsilon I)\#(Y+\varepsilon I).
\]
We use two standard properties of the geometric mean. If \(Z\ge0\)
and
\[
\begin{pmatrix}X&Z\\Z&Y\end{pmatrix}\ge0,
\]
then \(Z\le X\#Y\). Also,
\(X\#Y=X^{1/2}QY^{1/2}\) for some unitary \(Q\).

For \(A\in M_n(\mathbb C)\), define
\[
\D(A)=
\left\{
P\ge0:
\begin{pmatrix}P&A\\A^*&P\end{pmatrix}\ge0
\right\}.
\]

\begin{lemma} \label{lem:positive-blocks}
Let $A \in \M$ and let \(A=U|A|\) be a polar decomposition with \(U\) unitary.
\begin{enumerate}
\item
If \(X,Y\in\M\) are positive semidefinite and satisfy
\(|A^*|\le X\) and \(|A|\le Y\), then
\[
\begin{pmatrix}
X & A \\ A^* & Y
\end{pmatrix} \geq 0.
\]
In particular, if \(P\ge |A|\) and \(P\ge |A^*|\), then \(P\in\D(A)\).
\item
If \(P\in\D(A)\), then
\[
|A|\le P\#(U^*PU)\le \frac{P+U^*PU}{2}.
\]
Moreover, \(s(A)\prec_{w\log}s(P)\). In particular,
\(\|A\|\le\|P\|\) for every unitarily invariant norm. For
\(1\le i,j\le n\) with \(i+j-1\le n\), one also has
\[
s_{i+j-1}(A)^2\le s_i(P)s_j(P).
\]
\item
If \(P_k\in\D(A_k)\) for \(1\le k\le m\), then
\[
\sum_{k=1}^m P_k\in\D\!\left(\sum_{k=1}^m A_k\right).
\]
\end{enumerate}
 \end{lemma}

\begin{proof}
Let $A= U|A|$ be the polar decomposition of $A$.
For the first assertion, use the factorization
\[
\begin{pmatrix}  |A^*| & A  \\  A^*& |A|  \end{pmatrix} =
\begin{pmatrix}  U   \\  I \end{pmatrix}
|A|
\begin{pmatrix}  U^* & I \end{pmatrix} \ge 0.
\]

Now
\[
\begin{pmatrix}
X & A \\ A^* & Y
\end{pmatrix}
=\begin{pmatrix}
X-|A^*| & 0 \\ 0 & Y-|A|
\end{pmatrix}
+ \begin{pmatrix}
|A^*| & A \\ A^* & |A|
\end{pmatrix}
\geq 0.
\]

For the second assertion, conjugation by
\(\operatorname{diag}(U,I)\) gives
\[
\begin{pmatrix}U^*PU&|A|\\|A|&P\end{pmatrix}\ge0.
\]
The maximal characterization of the geometric mean and the
arithmetic--geometric mean inequality give the two Loewner bounds.
Set \(H=P\#(U^*PU)\). There is a unitary \(Q\) such that
\[
H=P^{1/2}Q(U^*PU)^{1/2}.
\]
Since \(|A|\le H\), eigenvalue monotonicity gives the first
inequality below. The multiplicative singular value inequality
applied to the factorization of \(H\) gives the second.
\[
\prod_{j=1}^k s_j(A)
\le \prod_{j=1}^k s_j(H)
\le \prod_{j=1}^k s_j(P),
\qquad 1\le k\le n.
\]
Weak log-majorization implies weak majorization and hence the norm
inequality. The bound for \(s_{i+j-1}(A)\) follows from
\(s_{i+j-1}(A)\le s_{i+j-1}(H)\) and the product inequality
\(s_{i+j-1}(XY)\le s_i(X)s_j(Y)\) applied to the same
factorization of \(H\).

The third assertion follows by summing the positive blocks
associated with \(P_k\in\D(A_k)\).
\end{proof}

\section{Main Results}

For \(1\le p\le\infty\), write \(R_p(A)=2^{1/p}|A|_p\).
For finite \(p\), \(R_p(A)=(|A|^p+|A^*|^p)^{1/p}\), and
\(R_\infty(A)=|A|_\vee\).

\begin{lemma}\label{L1}
For \(A\in M_n(\mathbb C)\) and \(1\le p\le\infty\), one has
\(2^{1/p}|A|_p\in\D(A)\).
\end{lemma}

\begin{proof}
Note that for finite $p$,
\[
|A|^p   \le |A|^p+|A^*|^p,
\qquad
|A^*|^p \le |A|^p+|A^*|^p .
\]
Hence, operator monotonicity of \(t^{1/p}\) gives
\[
|A|\le (|A|^p+|A^*|^p)^{1/p},
\qquad
|A^*|\le (|A|^p+|A^*|^p)^{1/p}.
\]
Since \((|A|^p+|A^*|^p)^{1/p}=2^{1/p}|A|_p\),
Lemma~\ref{lem:positive-blocks} gives
\[
\begin{pmatrix}
2^{1/p}|A|_p & A\\
A^* & 2^{1/p}|A|_p
\end{pmatrix}\ge0.
\]
As \(p\to\infty\), the diagonal blocks converge in norm to
\(|A|_\vee\). The limit block matrix is positive semidefinite,
which proves the assertion at \(p=\infty\).
\end{proof}

\begin{theorem}\label{M-C}
If \(A=U|A|\) with \(U\) unitary, then for \(1\le p\le\infty\),
\[
|A|\le 2^{1/p}\bigl(|A|_p\#(U^*|A|_pU)\bigr)
\le \frac{2^{1/p}}{2}\bigl(|A|_p+U^*|A|_pU\bigr).
\]
\end{theorem}

\begin{proof}
By Lemma~\ref{L1}, \(R_p(A)\in\D(A)\). Apply the second part of
Lemma~\ref{lem:positive-blocks} with \(P=R_p(A)\), then use joint
homogeneity of the geometric mean.
\end{proof}

\begin{theorem}\label{sums}
Let \(A_1,\ldots,A_m\in M_n(\mathbb C)\) and let \(V\) be a
unitary polar factor of \(\sum_{k=1}^m A_k\). For \(1\le p\le\infty\),
\begin{equation}\label{geom}
\begin{aligned}
\left|\sum_{k=1}^m A_k\right|&\le 2^{1/p}\biggl(\Bigl(\sum_{k=1}^m |A_k|_{p}\Bigr)\#\Bigl(\sum_{k=1}^m V^*|A_k|_{p}V\Bigr)\biggr)\\
&\le \frac{2^{1/p}}{2}\sum_{k=1}^m\bigl(|A_k|_p+V^*|A_k|_pV\bigr).
\end{aligned}
\end{equation}
The coefficient \(2^{1/p}\) in the geometric-mean bound is sharp
in each dimension \(n\ge2\).
\end{theorem}

\begin{proof}
Lemma~\ref{L1} gives \(R_p(A_k)\in\D(A_k)\) for every \(k\).
Adding the corresponding positive blocks gives
\[
2^{1/p}\sum_{k=1}^m |A_k|_{p}\in\D\!\left(\sum_{k=1}^m A_k\right).
\]
Inequality \eqref{geom} follows from Lemma~\ref{lem:positive-blocks}
and joint homogeneity of the geometric mean.
To prove sharpness, take \(m=1\) and, in dimension three, let
\[
A=\begin{pmatrix}
0&1&0\\
0&0&0\\
0&0&0
\end{pmatrix},
\qquad
V=\begin{pmatrix}
0&1&0\\
0&0&1\\
1&0&0
\end{pmatrix}.
\]
Then \(V\) is unitary and \(A=V|A|\), where
\(|A|=\operatorname{diag}(0,1,0)\). Since
\(|A^*|=\operatorname{diag}(1,0,0)\), we have
\[
|A|_p=2^{-1/p}\operatorname{diag}(1,1,0),
\qquad
V^*|A|_pV=2^{-1/p}\operatorname{diag}(0,1,1).
\]
The two matrices are diagonal, so their geometric mean is
computed entrywise. Hence
\[
2^{1/p}\bigl(|A|_p\#(V^*|A|_pV)\bigr)
=\operatorname{diag}(0,1,0)=|A|.
\]
This gives equality in the first inequality of \eqref{geom},
so the coefficient \(2^{1/p}\) is sharp. For \(n>3\), use
\(A\oplus 0_{n-3}\) and \(V\oplus I_{n-3}\).

For \(n=2\), take \(A=E_{12}\). Then \(|A|_p=2^{-1/p}I_2\),
and evaluating the bound on the second coordinate vector again
forces the coefficient to be at least \(2^{1/p}\).
\end{proof}

\begin{corollary}\label{cor:singular-values}
Under the hypotheses of Theorem~\ref{sums}, for \(i,j\ge1\) with
\(i+j-1\le n\),
\[
s_{i+j-1}\!\left(\sum_{k=1}^m A_k\right)
\le 2^{1/p}\sqrt{
s_i\!\left(\sum_{k=1}^m |A_k|_{p}\right)s_j\!\left(\sum_{k=1}^m |A_k|_{p}\right)}.
\]
The geometric mean of the two singular values is bounded above
by their arithmetic mean. In particular,
\[
s_{2j-1}\!\left(\sum_{k=1}^m A_k\right)
\le 2^{1/p}s_j\!\left(\sum_{k=1}^m |A_k|_{p}\right),
\qquad 2j-1\le n.
\]
Moreover,
\[
s\!\left(\sum_{k=1}^m A_k\right)
\prec_{w\log}2^{1/p}s\!\left(\sum_{k=1}^m |A_k|_{p}\right),
\]
and the same weak majorization holds with \(\prec_w\).
Consequently, for every unitarily invariant norm,
\begin{equation}\label{normsum}
\left\|\sum_{k=1}^m A_k\right\|\le 2^{1/p}\Bigl\|\sum_{k=1}^m |A_k|_{p}\Bigr\|.
\end{equation}
In each dimension \(n\ge2\), the coefficient in \eqref{normsum}
is sharp as a constant uniform over all unitarily invariant norms.
\end{corollary}

\begin{proof}
Apply Lemma~\ref{lem:positive-blocks} to
\(2^{1/p}\sum_{k=1}^m|A_k|_p\in\D(\sum_{k=1}^m A_k)\).
The odd-index estimate follows by setting \(i=j\). Weak
log-majorization implies weak majorization and hence \eqref{normsum}.
For \(m=1\) and \(A=E_{12}\in M_2(\mathbb C)\), one has
\[
\|E_{12}\|_{\mathrm{op}}
=1=2^{1/p}\bigl\||E_{12}|_p\bigr\|_{\mathrm{op}}.
\]
Thus equality holds in \eqref{normsum} for the operator norm.
For \(n>2\), the same conclusion follows by taking
\(A=E_{12}\oplus 0_{n-2}\).
\end{proof}

\subsection{Norm comparisons for a single matrix}

\begin{lemma}\label{B=Lemma}
Let \(X,Y\ge0\), and let \(\|\cdot\|\) be a unitarily invariant
norm. For a nonnegative concave function \(f\) on \([0,\infty)\),
\[
\|f(X+Y)\|\le \|f(X)+f(Y)\|.
\]
For a nonnegative convex function \(f\) on \([0,\infty)\) with
\(f(0)=0\), the inequality is reversed.
\end{lemma}

The concave case is due to Bourin and Uchiyama \cite{Bourin}.
For the convex case, see \cite[Corollary~2.2]{BL-concave}.

\begin{theorem}\label{N-Mo}
Let \(A\in M_n(\mathbb C)\). The map \(p\mapsto |A|_p\) is increasing
in the Loewner order on \([1,\infty]\). In particular,
\[
|A|_{\mathrm{sym}}\le |A|_p\le |A|_{\mathrm{qsym}}
\quad (1\le p\le2),\qquad
|A|_{\mathrm{sym}}\le |A|_{\mathrm{qsym}}\le |A|_p
\quad (2\le p\le\infty).
\]
For every unitarily invariant norm,
\[
\bigl\||A|_p\bigr\|
\le 2^{1-1/p}\bigl\||A|_{\mathrm{sym}}\bigr\|
\quad (1\le p\le\infty),
\qquad
\bigl\||A|_p\bigr\|
\ge 2^{1-1/p}\bigl\||A|_{\mathrm{sym}}\bigr\|
\quad (0<p\le1).
\]
Both factors are sharp for each unitarily invariant norm in
dimensions \(n\ge2\).
\end{theorem}

\begin{proof}
For \(1\le p\le q<\infty\), operator concavity of \(t^{p/q}\) gives
\[
\frac{|A|^p+|A^*|^p}{2}
\le \left(\frac{|A|^q+|A^*|^q}{2}\right)^{p/q}.
\]
Applying the operator-monotone power \(t^{1/p}\) proves monotonicity.
The endpoint \(q=\infty\) follows by taking limits.

Apply Lemma~\ref{B=Lemma} with \(X=|A|^p\), \(Y=|A^*|^p\), and
\(f(t)=t^{1/p}\). For \(p\ge1\), concavity gives
\(\|R_p(A)\|\le\||A|+|A^*|\|\), with the endpoint \(p=\infty\)
following by taking the norm limit. For \(0<p\le1\), convexity gives
the reverse comparison of \(\|(|A|^p+|A^*|^p)^{1/p}\|\) and
\(\||A|+|A^*|\|\). Dividing by \(2^{1/p}\) gives the stated bounds.
For \(A=E_{12}\) one has \(|A|_p=2^{-1/p}I_2\) and
\(|A|_{\mathrm{sym}}=I_2/2\), which gives equality in both estimates
for every unitarily invariant norm. For \(n>2\), the same argument applies to \(A\oplus 0_{n-2}\),
where \(A\) is the \(2\times2\) matrix above.
\end{proof}

\begin{corollary}\label{cor:power-norm-chains}
For every \(A\in M_n(\mathbb C)\) and every unitarily invariant norm,
\begin{align*}
\bigl\||A|_{\mathrm{sym}}\bigr\|
&\le \bigl\||A|_p\bigr\|
\le \bigl\||A|_{\mathrm{qsym}}\bigr\|
\le \sqrt{2}\,\bigl\||A|_{\mathrm{sym}}\bigr\|,
&&1\le p\le2,\\
\bigl\||A|_{\mathrm{sym}}\bigr\|
&\le \bigl\||A|_{\mathrm{qsym}}\bigr\|
\le \bigl\||A|_p\bigr\|
\le 2^{1-1/p}\bigl\||A|_{\mathrm{sym}}\bigr\|,
&&2\le p\le\infty.
\end{align*}
\end{corollary}

\begin{proof}
Use the Loewner comparisons from Theorem~\ref{N-Mo} and its upper
norm bound, taking \(p=2\) for the last term in the first chain.
\end{proof}

\subsection{Comparison with the quadratic bound}

The operator and norm estimates in Theorem~\ref{sums} and
Corollary~\ref{cor:singular-values} depend on the scaled sum
\[
2^{1/p}\sum_{k=1}^m |A_k|_p=\sum_{k=1}^m R_p(A_k).
\]
Although \(\sum_{k=1}^m|A_k|_p\) increases with \(p\) in the Loewner order,
the factor \(2^{1/p}\) decreases. A comparison with the
quadratic bound therefore requires control of the scaled
sum. Such a comparison follows from the scalar power
inequality when each pair \(|A_k|\), \(|A_k^*|\) commutes.

A matrix \(A\in M_n(\mathbb C)\) is called \emph{binormal} if
\(A^*A\) and \(AA^*\) commute. Equivalently, \(|A|\) and
\(|A^*|\) commute.

\begin{example}
\begin{enumerate}
\item
Let \(A=PD\), where \(P\) is a permutation matrix and \(D\ge0\)
is diagonal. Then
\[
|A|=D
\qquad\text{and}\qquad
|A^*|=PDP^*.
\]
Both moduli are diagonal, so \(A\) is binormal. It is nonnormal
whenever \(PDP^*\ne D\).
\item
Let \(A\) be a square-zero matrix, so that \(A^2=0\). Then
\[
\operatorname{ran}|A^*|
=\operatorname{ran}A
\subseteq\ker A
=\ker|A|.
\]
Thus \(|A||A^*|=0\), and taking adjoints also gives
\(|A^*||A|=0\). Hence \(A\) is binormal. Since the two moduli
have orthogonal supports,
\[
R_p(A)=|A|+|A^*|,
\qquad 1\le p\le\infty.
\]
\end{enumerate}
\end{example}

\begin{proposition}\label{prop:commuting-moduli}
Let \(A_1,\ldots,A_m\in M_n(\mathbb C)\) be binormal matrices.
Then the scaled sum \(2^{1/p}\sum_{k=1}^m |A_k|_p\) is
nonincreasing in the Loewner order. More precisely,
\[
2^{1/q}\sum_{k=1}^m |A_k|_q
\le 2^{1/p}\sum_{k=1}^m |A_k|_p,
\qquad 1\le p\le q\le\infty.
\]
Consequently, if \(V\) is a unitary polar factor of
\(\sum_{k=1}^m A_k\), then for \(2\le q\le\infty\),
\[
\begin{aligned}
&2^{1/q}
 \Bigl(\sum_{k=1}^m |A_k|_q\Bigr)
 \#\Bigl(\sum_{k=1}^m V^*|A_k|_qV\Bigr)\\
&\qquad\le
 \sqrt{2}
 \Bigl(\sum_{k=1}^m |A_k|_2\Bigr)
 \#\Bigl(\sum_{k=1}^m V^*|A_k|_2V\Bigr).
\end{aligned}
\]
In particular,
\[
2^{1/q}\Bigl\|\sum_{k=1}^m |A_k|_q\Bigr\|
\le \sqrt{2}\Bigl\|\sum_{k=1}^m |A_k|_2\Bigr\|
\]
for every unitarily invariant norm.
\end{proposition}

\begin{proof}
Fix \(k\). Since \(A_k\) is binormal, the matrices \(|A_k|\)
and \(|A_k^*|\) are simultaneously diagonalizable. In a common
orthonormal eigenbasis, each diagonal entry of \(R_p(A_k)\)
has the form \((a^p+b^p)^{1/p}\), where \(a,b\ge0\).
For \(p\le q<\infty\), put \(d=(a^p+b^p)^{1/p}\). If \(d>0\),
then \(a/d,b/d\in[0,1]\), so
\[
(a/d)^q+(b/d)^q\le (a/d)^p+(b/d)^p=1.
\]
Thus \((a^q+b^q)^{1/q}\le d\). The case \(d=0\) is immediate,
and the case \(q=\infty\) follows by taking the limit.
Therefore,
\[
2^{1/q}|A_k|_q\le 2^{1/p}|A_k|_p,
\qquad 1\le p\le q\le\infty.
\]
Summing over \(k\) proves the first inequality. Taking \(p=2\)
and using joint homogeneity and monotonicity of the geometric
mean gives the second inequality. The norm comparison follows
from monotonicity of unitarily invariant norms on positive
semidefinite matrices.
\end{proof}

Proposition~\ref{prop:commuting-moduli} shows that the bounds are
nonincreasing in the exponent for binormal summands. The following
example shows that without binormality, both bounds can increase
from \(p=2\) to \(p=4\), even for two \(2\times2\) summands.

\begin{example}
Let
\[
A_1=\frac14
\begin{pmatrix}
3&-\sqrt3\\
\sqrt3&-1
\end{pmatrix},
\qquad
A_2=\frac14
\begin{pmatrix}
-1&\sqrt3\\
-\sqrt3&3
\end{pmatrix}.
\]
The matrices \(A_1^*A_1\) and \(A_1A_1^*\) are rank-one
projections. In fact,
\[
A_1^*A_1=\frac14\begin{pmatrix}3&-\sqrt3\\-\sqrt3&1\end{pmatrix},
\qquad
A_1A_1^*=\frac14\begin{pmatrix}3&\sqrt3\\\sqrt3&1\end{pmatrix}.
\]
Each is positive semidefinite with trace one and determinant zero.
A projection is unchanged by taking any positive power, so
\(|A_1|^p+|A_1^*|^p=\operatorname{diag}(3/2,1/2)\). Hence
\[
R_p(A_1)=\operatorname{diag}\bigl((3/2)^{1/p},(1/2)^{1/p}\bigr).
\]
Since \(A_2\) is obtained from \(A_1\) by interchanging the
two coordinates, the diagonal entries of \(R_p(A_2)\) are
the same in reverse order. Therefore,
\[
R_p(A_1)+R_p(A_2)
=\bigl((3/2)^{1/p}+(1/2)^{1/p}\bigr)I_2.
\]
The scalar coefficient is less than \(1.94\) at \(p=2\)
and greater than \(1.94\) at \(p=4\). Indeed,
\[
\sqrt{3/2}<1.23,\qquad \sqrt{1/2}<0.71,
\]
whereas
\[
(3/2)^{1/4}>1.10,\qquad (1/2)^{1/4}>0.84.
\]
Thus the scaled sum increases strictly in the Loewner order
from \(p=2\) to \(p=4\). Since both scaled sums are scalar
multiples of \(I_2\), the geometric-mean and unitarily invariant
norm bounds also increase strictly.
\end{example}

\begin{example}\label{ex:strict-power-comparison}
Let
\[
A=\begin{pmatrix}0&2\\1&0\end{pmatrix}.
\]
Then \(|A|=\operatorname{diag}(1,2)\) and
\(|A^*|=\operatorname{diag}(2,1)\). Thus \(A\) is nonnormal,
its two moduli commute, and
\[
R_p(A)=(1+2^p)^{1/p}I_2,\qquad R_\infty(A)=2I_2.
\]
For every finite \(p>2\), the inequality
\((1+2^p)^{1/p}<\sqrt{5}\) gives a strictly smaller
geometric-mean upper bound than the quadratic estimate.
At \(p=\infty\), the operator-norm bound is exact.
\end{example}

\section*{Declaration of generative AI and AI-assisted technologies in the writing process}

During the preparation and revision of this manuscript, the authors used
ChatGPT for language editing, \LaTeX{} formatting, manuscript organization,
consistency checks, and assistance with computations and numerical
verification. The authors retained full control over the mathematical
development, analysis, and presentation of the results and take full
responsibility for the content of the publication.

\bibliographystyle{plain}

\end{document}